\documentclass[11pt]{amsart}
\usepackage{amsmath,amssymb}

\theoremstyle{plain}
\newtheorem{propn}{Proposition}[section]
\newtheorem{thm}[propn]{Theorem}

\theoremstyle{definition}
\newtheorem{defn}[propn]{Definition}

\theoremstyle{remark}

\numberwithin{equation}{section}

\begin{document}

\title[Real hypersurfaces of hyperbolic spaces]{Real hypersurfaces of complex and quaternionic
hyperbolic spaces}

\author{Thomas Murphy}
\address{School of Mathematical Sciences, University College Cork, Ireland.}

\curraddr{D\'{e}partment de Math\'{e}matique,
Universit\'{e} Libre de Bruxelles,
 \ Boulevard du Triomphe,
B-1050 Bruxelles,
Belgique.}

\email{tmurphy@ulb.ac.be}


\date{April 1st, 2012.}


\maketitle

\thispagestyle{empty}

\begin{abstract}
We introduce curvature-adapted foliations of complex hyperbolic space and study some of their properties. Generalized
pseudo-Einstein hypersurfaces of complex hyperbolic space are classified. Analogous results for curvature-adapted
foliations of quaternionic hyperbolic space are also discussed.
\end{abstract}

\section{Introduction}

Curvature-adapted hypersurfaces are of great interest to geometers and have been the focus of much attention since the
concept was introduced by d'Atri \cite{datri}. Into this class fall many known examples of hypersurfaces with constant
principal curvatures, a subject intensively studied since the time of Levi-Civita and Segre. Generalizing many distinguished
families of hypersurfaces (for example umbilic hypersurfaces), their geometry is adapted to that of the ambient space in a
special way. The Riccati equation was used by Gray in his monograph \cite{gray} to develop foundational theorems about
curvature-adapted submanifolds in complex space forms: here it is employed to tackle outstanding problems in complex and
quaternionic hyperbolic spaces.

Throughout this paper, we define $M$ to be a connected hypersurface of a connected, simply connected rank one symmetric space $\overline{M}$, $\overline{R}$ to be the Riemannian curvature
tensor of $\overline{M}$ and $\xi$ to be a unit normal vector of $M$ at $p\in M$. We give the Riemannian metrics of $\overline{M}$ the standard scaling, so that  sectional curvatures lie between $\pm1$ and $\pm4$. The normal Jacobi operator
$$
K_{\xi} := \overline{R}(\xi, \cdot)\xi \in End(T_pM)
$$
of $M$ (with respect to $\xi$) describes the curvature of the ambient manifold $\overline{M}$ at $p$, whereas the shape
operator $A_{\xi}$ of $M$ (with respect to $\xi$) describes the curvature of $M$ as a submanifold of $\overline{M}$ in
direction $\xi$. Both of these are self-adjoint operators, and hence have eigendecompositions. $M$ is said to be curvature-adapted if these operators are simultaneously diagonalizable, that is if
$$
K_{\xi}\circ A_{\xi} = A_{\xi}\circ K_{\xi}
$$
at every point $p\in M$. This means that a common eigenbasis for $K_{\xi}$ and $A_{\xi}$ exists at every point, which will
generically be denoted $E$.

The geometry of curvature-adapted hypersurfaces in rank one symmetric spaces has been a fruitful field of study and there is
a substantial body of literature concerned with their classification. In real space forms it is easy to see that every
hypersurface is curvature-adapted. In non-flat complex space forms they coincide exactly with Hopf hypersurfaces. Examples abound; any
tube around any complex submanifold of complex projective space is a Hopf hypersurface \cite{cecilryan}. Recall that a Hopf
hypersurface is  defined by the property that the structure vector field $-J\xi$ is a principal curvature vector field. For
a Hopf hypersurface with structure vector field $U = -J\xi$, denote the principal curvature function corresponding to the
structure vector field by $\alpha$. It is well-known that $\alpha$ is constant for $\mathbb{C}P^n$ \cite{alpha1} and $\mathbb{C}H^n$ \cite{alpha2}.  An explicit classification for Hopf hypersurfaces has been achieved in $\mathbb{C}P^n$
\cite{kimura} and $\mathbb{C}H^n$ \cite{berndt1} under the assumption of constant principal curvatures.

\begin{defn}
A singular Riemannian foliation $\mathcal{F}$ of $\overline{M}$ is a decomposition of $\overline{M}$  into 
embedded submanifolds, called the leaves of the foliation $L$, such that
\begin{enumerate}
\item $ T_pL = \lbrace X_p, X \in \Xi_{\mathcal{F}}\rbrace $ for every $L\in \mathcal{F}$ and $p\in L$, where $\Xi_{\mathcal{F}}$
is the space of smooth vector fields on $\overline{M}$ everywhere tangent to the leaves of $\mathcal{F}$, and
\item every geodesic is orthogonal to the leaves at all or none of its points.
\end{enumerate}
\end{defn}

The leaves of maximal dimension are called regular, otherwise they are singular. If every leaf is regular one recovers the
traditional definition of a Riemannian foliation. Let $\nu_pL$ denote the normal space to a leaf $L$ at a point $p$.  A Riemannian foliation is said to admit sections (or to be \emph{polar}) if there 
exists a complete immersed submanifold  $\Sigma\subset \overline{M}$ such that, for all leaves $L\in \mathcal{F}$, $L\cap \Sigma \neq \emptyset$ and,
for all points $p\in \Sigma$, one has $T_p\Sigma \subset \nu_pL$. Singular Riemannian foliations arise naturally in geometry; for example in the
study of isometric group actions and Riemannian submersions.

\begin{defn}
We define a connected submanifold $P\subset \overline{M}$ to be curvature-adapted if, for every unit normal vector $\xi$ at a point
$p$, 
\begin{enumerate}
\item  $\overline{R}(\xi, X)\xi\in T_pP$ for $X\in T_pP$, and \\
\item  if $K_{\xi} := \overline{R}|_{T_pP}$, then $K_{\xi}$ and $A_{\xi}$ commute, 
\end{enumerate}
at all points $p\in P$. 
\end{defn}

This generalizes the definition
for hypersurfaces. Let $\nu^1(P)$ be the unit sphere bundle. For $t\in \mathbb{R}^+$, set $M_t:= \lbrace exp(t\xi)\text{ : } \xi\in\nu^1(P)\rbrace$. This is the tube of radius $t$ around $P$.  Gray's Theorem (Theorem (\ref{gt}) here) states that the tubes of sufficiently small radius around a curvature-adapted submanifold
$P\subset \overline{M}$ are curvature-adapted hypersurfaces. Motivated by this result,  let us make the following
definition:
\begin{defn}
A \emph{curvature-adapted foliation} $\mathcal{F}$ of $\overline{M}$ is a singular Riemannian foliation of $\overline{M}$
whose regular leaves arise as the tubes around a curvature-adapted submanifold $P\in \mathcal{F}$ (together with $P$ if there are no singular leaves).
\end{defn}

The horospherical foliation of $\mathbb{C}H^n$ is an example of a curvature-adapted foliation without singular leaves. Recall, following Ivey and Ryan \cite{ivey}, \cite{iveyryan}, a Hopf hypersurface of $\mathbb{C}H^n$ is said to
have \emph{small Hopf curvature} if $0< \alpha < 2$. They have constructed many examples of such hypersurfaces.
Similarly one defines a Hopf hypersurface of $\mathbb{C}H^n$ to have large Hopf curvature if $2 < \alpha$. Any
tube around a complex submanifold gives a Hopf hypersurface with large Hopf curvature. In the borderline case ( i.e.
$\alpha^2 - 4 = 0$ ) one says that the hypersurface is degenerate.

A curvature-adapted foliation $\mathcal{F}$  of $\mathbb{C}H^n$ is said to be degenerate (resp.
non-degenerate) if a regular leaf $M\in \mathcal{F}$ is degenerate (resp. non-degenerate). It follows from the Riccati equation that
if $M$ is a degenerate (resp. non-degenerate) Hopf hypersurface, then so are all parallel hypersurfaces $M_t$.

We shall show the following:
\begin{thm}\label{thm2}
A regular leaf of a curvature-adapted foliation $\mathcal{F}$ of  $\mathbb{C}H^n$ has small Hopf curvature if, and only if,
$\mathcal{F}$ is the set of all tubes around a totally geodesic $\mathbb{R}H^n$ together with $\mathbb{R}H^n$. 
\end{thm}

A precisely analogous result also holds true for curvature-adapted foliations of $\mathbb{H}H^n$, as there are
similar results on the spectral data of the shape operator known \cite{berndt}.  Some partial results for degenerate
curvature-adapted foliations of $\mathbb{C}H^n$ and $\mathbb{H}H^n$ are then given.

We  present in Section $5$ new proofs of the work of Okumura \cite{okumura}, and Montiel-Romero \cite{mr}, classifying real
hypersurfaces of $\mathbb{C}P^n$ and $\mathbb{C}H^n$ whose induced almost-contact   structure $P$ commutes with $A_{\xi}$. Our approach,
exploiting the spectral data of the shape operator, is considerably shorter and also works for the analogous problem in
quaternionic space forms $\mathbb{H}H^n$. This unifies their work with that of Lyu-Perez-Suh \cite{lyuperezsuh}.

To conclude we investigate a related problem.  It is well-known that an Einstein manifold cannot be embedded isometrically
as a hypersurface in complex projective space $\mathbb{C}P^n$. Kon \cite{kon} defined a hypersurface to be
\emph{pseudo-Einstein} if there exist constants $\rho$ and $\sigma$ so that for any tangent vector $X$,
$$
SX = \rho X + \sigma\langle X,U\rangle U
$$
where $S$ denotes the (1,1)-Ricci tensor and $U = -J\xi$ denotes the structure vector field. Kon classified such
hypersurfaces under the assumption $n\geq 3$. Montiel \cite{montiel} derived a classification in complex hyperbolic space
$\mathbb{C}H^n$ under the same assumption. A canonical generalization of this concept is to allow $\rho$ and $\sigma$ to be
non-constant smooth functions. Such hypersurfaces are called generalized pseudo-Einstein. Cecil and Ryan \cite{cecilryan}
showed that in $\mathbb{C}P^n, n\geq 3$ such hypersurfaces coincide precisely with the pseudo-Einstein hyersurfaces, i.e.
$\rho$ and $\sigma$ must in fact be constant. In the survey paper of Niebergall and Ryan \cite{nieryan} the following open
problem is listed:
\begin{itemize}
\item The Cecil-Ryan theorem shows us that the assumption that $\sigma$
and $\rho$ are constant is unnecessary in Kon's work. Is the analogous statement true for complex hyperbolic space?
\end{itemize}
In other words, if $M\subset \mathbb{C}H^n, n\geq 3$ is a generalized pseudo-Einstein hyersurface, must it in fact be a
pseudo-Einstein hypersurface? If this were true, Montiel's work would classify such hypersurfaces in $\mathbb{C}H^n$, $n\geq
3$. Recently Kim and Ryan \cite{kimryan} have classified generalized pseudo-Einstein hypersurfaces in $\mathbb{C}P^2$ and
Ivey and Ryan \cite{iveyryan} classified such hypersurfaces in $\mathbb{C}H^2$, so the classification of generalized
pseudo-Einstein hypersurfaces in complex space forms with $n=2$ is settled. We establish:
\begin{thm}\label{cor2}
Let $M$ be a generalized pseudo-Einstein hypersurface of complex hyperbolic space $\mathbb{C}H^n, n\geq 3$. Then $M$ is
congruent to an open part of:
\begin{enumerate}
\item a tube of radius $r\in \mathbb{R}^+$ around a
totally geodesic $\mathbb{C}H^k\subset \mathbb{C}H^n$, where $k=0$ or $n-1$, or \\
\item a horosphere.
\end{enumerate}
\end{thm}

\section{Curvature-adapted foliations}

For  this section, we set our ambient manifold to be $\overline{M} = \mathbb{C}H^n$. We  denote the K\"ahler structure by $J$ and the Levi-Civita connection by $\overline{\nabla}$. The curvature tensor
$\overline{R}$ is given as
\begin{align*}
\overline{R}(X,Y)Z =  & -\langle Y,Z\rangle X + \langle X,Z\rangle Y\\
& - \langle JY,Z\rangle JX + \langle JX,Z\rangle JY +2 \langle JX,Y\rangle JZ 
\end{align*}

We first recall a theorem of Gray (\cite{gray}, Theorem 6.14):
\begin{thm}\label{gt}
Let $P\subset \overline{M}$ be an embedded curvature-adapted submanifold. Then
\begin{itemize}
\item any tube around $P$ is also curvature-adapted, and
\item the common eigenbasis $E(C_{\xi}(t))$ of the shape operator $A_{\xi}(t)$ of the tube around $P$ and the normal Jacobi operator $K_{\xi}(t)$  may be chosen parallel along geodesics normal to $P$.
\end{itemize}
\end{thm}

Given a hypersurface $M = M_0\subset \overline{M}$, the parallel hypersurface $M_a$ is defined by, for each
point $p\in M$, flowing by a distance $a$ along the geodesic $C_{\xi}(t) = exp_p(t\xi)$ which passes through $p$ with initial
direction $\xi$. The Riccati equation describes the evolution of the shape operator $A_{\xi}(t)$ of nearby parallel
hypersurfaces in terms of its derivative and the normal Jacobi operator $K_\xi$ along the geodesic $C_{\xi}(t)$:
$p\in M$:
\begin{equation}\label{Riccati}
A_{\xi}'(t) = (A_{\xi}(t))^2 + K_{\xi}(t),
\end{equation}
with $A_{\xi}(0)$ the shape operator of $M_0$ at $p$.
This equation relates the shape operator of nearby parallel hypersurfaces to $K_{\xi}$. Taking a curvature-adapted
hypersurface simplifies this to a family of ordinary differential equations, namely
$$\lambda_i'(t) = (\lambda_i(t))^2+ \kappa_i,
$$
for $i=1,\dots, dim(\overline{M})-1$, with initial conditions $\lambda_i(0)$ being the principal curvatures of $M_0$ at $p$. Here $\kappa_i = -1$ or $-2$.  This suggests that the investigation of nearby
parallel hypersurfaces to a given curvature-adapted hypersurface might be profitable. For the rest of this section we assume that $M$ is curvature-adapted.

Associated to $M$ is the induced Levi-Civita connection $\nabla$, and
we denote its curvature tensor by $R$. Let $TM$ denote the tangent bundle of $M$ and $\nu(M)$ the normal bundle.  Set $U=
-J\xi$ and let $\mathfrak{D}:= (\mathbb{R}U)^{\perp}$ denote the maximal complex subbundle of $TM$.

The Riccati equation for hypersurfaces is equivalent to the Jacobi equation, where there are also  well-established techniques to describe the principal curvatures of nearby parallel hypersurfaces. Recall that $J$ is said to be an $M$-Jacobi field if it is a non-zero Jacobi field along $C_{\xi}$ satisfying the initial conditions $J(0)\in T_pM$ and $J'(0) = -A_{\xi}(J(0))$. A focal point of $M$ along $C_{\xi}$ is given as $C_{\xi}(t_0)$ when $J(t_0)=0$ for an $M$-Jacobi vector field. We refer the reader to \cite{b}, Chapter 8 for further details about $M$-Jacobi vector fields. Take  $E_i(C_{\xi}(t)$ to be a basis for the $M$-Jacobi vector fields along $C_{\xi}$ from $M$.  As $M$ is curvature-adapted, $E_i$ being a focal point is equivalent to the corresponding principal curvature function developing a singularity (i.e. one has $\lambda_i(t_0) = \infty$). Let us now show how this is the case.

Suppose that $\mathcal{F}$ is a singular Riemannian foliation with a singular leaf $P$. We firstly explain how a regular leaf $M$ of  $\mathcal{F}$ and the singular leaf $P$ are related. 

 One can calculate the principal curvatures of the tube around $P$ using the  same calculation as in Lemma 7.8 of \cite{gray}.  The Riccati equation for the tube around a curvature-adapted submanifold $P$ simplifies to the following family of equations along $C_{\xi}(t)$, where $C_{\xi}(0)=p\in P$:
$$
\lambda_i'(t) = \lambda_i^2(t) + \kappa_i,
$$
for $i=1,\dots, dim(\overline{M})-1$,  $\kappa_i = -1$ or $-2$, and initial conditions
\begin{enumerate}
\item $\lambda_i(0) = \lambda_i(p)$, for $i=1, \dots, dim(P)$,
\item $\lambda_i(0) = -\infty$, for $i= dim(P)+1, \dots, 2n-1$.
\end{enumerate}

Now one can start with a regular leaf $M$ which we take to be the tube of radius $t_0$ around $P$. Then one can calculate the principal curvatures of the tube around $M$ along the same geodesic, but this time with initial conditions the point $q\in M$ and initial normal vector $\overline{\xi}$. Here $\overline{\xi}(0) = -\xi(t_0)$, etc.  Travelling back along $C_{\overline{\xi}}(t)$ to the point  $p\in P$, the vector fields in $E_i(t)$  whose  corresponding principal curvature functions ``focalize" at $P$ (i.e. the functions $\lambda_i(t)$  which become infinite at $p$) in the eigenbasis $E(C_{\overline{\xi}}(t))$ are precisely the $M$-Jacobi vector fields which have a focal point at the point $p$. This follows from the equation $E_i'(t) = -\lambda_i(t)(E_i(t)$.  We say $P$ is the focal leaf of $M$.

Let us state the following  well-known result:
\begin{thm}\cite{berndt1}
A hypersurface $M\subset \mathbb{C}H^n$ with constant principal curvatures is Hopf if and only if it is locally congruent to
an open part of;
\begin{itemize}
\item  a tube of radius $r\in \mathbb{R}^+$ around a totally
geodesic $\mathbb{C}H^{k}\subset \mathbb{C}H^{n}$ for $k = 0,\dots, n-1$, \\
\item  a tube of radius $r\in \mathbb{R}^+$ around a totally
geodesic $\mathbb{R}H^{n}\subset \mathbb{C}H^{n}$, \\
\item  a horosphere. \\
\end{itemize}
\end{thm}

Let  $\sigma_p(\mathfrak{D})$ denote the spectrum of
$A_{\xi}|_{\mathfrak{D}_p}$. Set $ T_{\lambda} = ker
(A_{\xi}|_{\mathfrak{D}_p} - \lambda id_p)$ to be the eigenspace
associated with $\lambda\in\sigma_p(\mathfrak{D})$.
\begin{thm} (See \label{key}\cite{berndt0}, \cite{boning}.)
\begin{enumerate}
\item Suppose $M\subset \mathbb{C}H^n$ is a degenerate Hopf hypersurface. Then \newline $1$
$\in \sigma_p(\mathfrak{D}), \text{ for all } p\in M$ and $JT_{\lambda}(p) \subset T_{1}(p)$  for all $p\in
M, \lambda \in \sigma_p(\mathfrak{D})\backslash \lbrace 1 \rbrace$.
\item Suppose $M\subset \mathbb{C}H^n$ is a  Hopf hypersurface with small Hopf curvature. Then given $\lambda \in
 \sigma_p(\mathfrak{D})$, an associated eigenvector $Y$
satisfies $A_{\xi}(JY) = \lambda^*(JY)$ where $\lambda^*$ satisfies the equation
\begin{equation}\label{bb}
(2\lambda^* - \alpha)(2\lambda -\alpha) = \alpha^2 -4.
\end{equation}
\end{enumerate}
\end{thm}

Choose $p\in M\subset \mathbb{C}H^n$, with
 corresponding section $C_{\xi}$ through a point $p$ and suppose $M$ is degenerate. Choose an
orthornormal framing $E_i(t)$ along $C_{\xi}(t)$ such that $E_1(t)= -J\xi(t)$, and for $E_{2k}(t)\in$ $T_{\lambda(t)}(t)$,
$\lambda(t) \neq 1$, $k\geq 1$, we take $E_{2k+1}(t) = JE_{2k}(t) \in T_1(t)$. We may always choose such
a framing by Theorem (\ref{key}), and this framing is parallel along $C_{\xi}$ (i.e. $\overline{\nabla}_{\xi}E_i = 0$).
Similarly if $M$ has small Hopf curvature choose $E_{2k}(t)\in T_{\lambda(t)}(t)$ and $E_{2k+1}(t) = JE_i(t)$ $\in T_{\lambda*(t)}(t)$, with $k\geq 1$.

In the above we have lightly abused notation by replacing $C_{\xi}(t)$ with $t$ for ease of exposition when the context is clear. We remark that for the case $\overline{M}= \mathbb{H}H^n$ results which are direct analogues of the above results are to found in \cite{berndt}. This allows one to prove precisely analogous results.

\section{Proof of Theorem \ref{thm2}}

 \proof The foliation induced by taking the  tubes of
radius $t>0$ around a totally geodesic $\mathbb{R}H^n\subset \mathbb{C}H^n$, together with the focal set $\mathbb{R}H^n$,
has small Hopf curvature. This may be calculated using the $M$-Jacobi field theory outlined, and is carried out in \cite{berndt1}.

 Suppose conversely $\mathcal{F}$ satisfies the assumptions of the theorem and
has a regular leaf $M$ with small Hopf curvature. Pick a point $p\in M$. Choose the framing $E_i(t)$ along $C_{\xi}(t)$ as explained in the last section. Then from Equation (\ref{bb}) we
see precisely half the principal curvatures lying in $\sigma_p(\mathfrak{D})$ must focalize. This is because small Hopf curvature implies that \newline $
(2\lambda_i - \alpha)(2\lambda_i^* - \alpha)  $
is now negative, so one of these two terms is positive. In particular, we may assume $\lambda_i > \frac{\alpha}{2}$ without loss of generality. Solving the corresponding Riccati equation yields
$$
\lambda_i(t) = Coth(\theta_i(p) -t),
$$
where $0 < \theta_i(p) < \infty$  is chosen so that $Coth(\theta_i)=\lambda_i(p)$. 
This focalizes at distance $t=\theta_i$. As $\mathcal{F}$ is a Riemannian foliation, the distance between this focal leaf and $M$ is constant so this implies that $\theta_i$ is independent of our choice of $p$. Thus $\lambda_i$ is locally constant and has constant multiplicity. Then Equation (\ref{bb}) implies $\lambda_i^*$ is also constant. Thus it follows that
$M$ is a Hopf hypersurface with constant principal curvatures and focal set a totally real submanifold, and we are done by the main result of \cite{berndt1}.
\endproof

\section{degenerate Hopf hypersurfaces}

The class of degenerate hypersurfaces merits further investigation as  the usual equations stemming from the
Gauss-Codazzi-Ricci equations break down. For this reason the usual techniques do not work, and so many classification
results for Hopf hypersurfaces assume non-degeneracy. B\"{o}ning \cite{boning} observes that the horosphere is the only
known degenerate Hopf hypersurface: the obvious question is if there are any more. B\"{o}ning \cite{boning} also classified
non-degenerate Hopf hypersurfaces with at most three principal curvatures in $\mathbb{C}^n, n\geq 3$; they are all
homogeneous. In recent work, Ivey and Ryan \cite{iveyryan} have associated a degenerate Hopf hypersurface $M^3$ of
$\mathbb{C}H^2$ to any contact curve in $\mathbb{S}^3$ using exterior differential systems. Their construction actually
associates a contact curve in $\mathbb{S}^3$ to an \emph{adapted lift} of $M$, $\tilde{M}\subset SU(2,1)$. In the degenerate
case $\tilde{M}$ depends on one function of one variable. However they cannot explicitly compute the principal curvatures of
$M^3$ from $\tilde{M}$ and so it is still an open problem to construct a degenerate Hopf hypersurface of $\mathbb{C}H^2$ (or
$\mathbb{C}H^n$ in general) which is not the horosphere.

It is obvious that the unique degenerate curvature-adapted foliation of $\mathbb{K}H^n$, $\mathbb{K} = \mathbb{C}$ or
$\mathbb{H}$ where the non-Hopf principal curvatures of the regular leaves satisfy  $|\lambda_i(p)| \geq 1$ is the horospherical
foliation, because from the Riccati equation one sees that each regular leaf must have constant principal curvatures and we have already remarked that such hypersurfaces are classified.  We also have the following result:
\begin{propn}
There is no degenerate curvature-adapted foliation of $\mathbb{C}H^n$  with a
unique minimal singular leaf.
\end{propn}

\proof  In our notation, $P$ denotes the minimal singular leaf and $M$ a regular leaf. Fix $p\in M$ and consider the section $C_{\xi}$
passing through $p$, where $\xi$ points in the direction of $P$. Let $\theta_1$ denote the distance between $M$ and $P$ along
$C_{\xi}$. As each regular leaf of $\mathcal{F}$ is degenerate, the Hopf principal curvature $\alpha=2$ is constant along $C_{\xi}$
and does not focalize at $P$. Let $\lambda_1 = Coth(\theta_1)\in \sigma_p(\mathcal{D})$ denote  the principal curvature of the shape operator of $M$ which focalizes at
$P$. Then by Theorem (\ref{key}) $JT_{\lambda_1}(0) \in T_{1}(0)$. This vector subspace is invariant under
parallel translation along $C_{\xi}$ and so corresponds to a subspace of $T_{C_{\xi}(\theta_1)}P$. Any other principal curvature
function $\lambda_i, i > 1$ of $X$ at $T_{C_{\xi}(\theta_1)}X$ has corresponding eigenspace $T_{\lambda_i}$, and from the known
spectral data $JT_{\lambda_i}\subset T_{1}$. But $|\lambda_i| < 1$ for $i > 1$ as there is
a unique singular leaf. Otherwise it would follow from the Riccati equation that there is another solution of the shape operator of the regular leaves of the form $\lambda_2(t)= Coth(\theta_2 -t)$, $\theta_2\neq \theta_1$ and so a second focal submanifold would exist.  From this it follows that $P$ cannot be minimal. To see this, observe that at least half of the principal curvatures in $\sigma(\mathfrak{D})$ are $+1$, so they outweigh the other principal curvatures which all have $|\lambda_i|<1$ and so $P$ cannot be minimal. \endproof

Again there is an analogous result for $\mathbb{H}H^n$.

\section{Applications to the study of almost-contact structures} 

In this section, $\overline{M}$ is either $\mathbb{C}P^n$ or $\mathbb{C}H^n$, and $M\subset \overline{M}$ is a hypersurface.
Let $P$ denote the
skew-symmetric $(1,1)$ tensor field on $M$ given by
$$
JX = PX + \langle X, U\rangle \xi.
$$
This is the induced \emph{almost contact structure} on $TM$. By the Gauss formula and the  Weingarten equation we obtain
$$
\nabla_X U = PA_{\xi}(X).
$$

Recall the well-known  fact that the shape operator of a real hypersurface of $\mathbb{C}H^n$ or $\mathbb{C}P^n$ cannot vanish: in fact
$$
\| \nabla A\|^2 \geq  4(n-1)
$$
in the standard scaling. Equality in the above bound is achieved precisely by hypersurfaces whose shape operators and induced almost-contact
structure commute;
$$
A_{\xi} \circ P = P \circ A_{\xi}.
$$

Thus classifying which hypersurfaces in non-flat complex space forms achieve equality is a natural question. We now outline a
simplified proof of the work of Okumura \cite{okumura} and Montiel and Romero \cite{mr}, who answered this question.

If $P \circ A_{\xi} = A_{\xi}\circ P$ then it is easy to see $M$ is Hopf. We need now the fact that, for $\overline{M}=\mathbb{C}P^n$ there is an exact analogue of the framing along $C_{\xi}(t)$ and, as shown in \cite{boning}, Equation (\ref{bb}) also holds for this framing. Suppose that $M$ is non-degenerate (the degenerate case is analogous, using Theorem (\ref{key})). Then given $\lambda \in \sigma_p(\mathfrak{D})$, Equation (\ref{bb}) together
with the equation  $P \circ A_{\xi} = A_{\xi}\circ P$ implies $\lambda=\lambda^*$. Thus
$$
(2\lambda - \frac{\alpha}{2})^2 = \alpha^2 \pm 4,
$$
so $\lambda$ is constant.  Hence $M$ is a Hopf hypersurface with constant principal curvatures of a non-flat complex space form, and a case-by-case check shows that $M$ is isometric to one of
\begin{enumerate}
\item a tube of radius $t$ around a totally geodesic $\mathbb{C}P^k\subset \mathbb{C}P^n$, $0\leq k \leq n-1$, $0< t <
\frac{\pi}{2}$,
\item a tube of radius $t$ around a totally geodesic $\mathbb{C}H^k\subset \mathbb{C}H^n$, $0\leq k \leq n-1$, $0< t <
\infty$,
\item a horosphere of $\mathbb{C}H^n$.
\end{enumerate}

We note that in \cite{lyuperezsuh} the authors investigated in
$\mathbb{H}H^n$ an analogous problem to that studied by Romero and Montiel in $\mathbb{C}H^n$. Again, every real hypersurface $M\subset \mathbb{H}H^n$ satisfies a similar bound.
They prove that equality in the bound is achieved precisely by hypersurfaces satisfying
$$A_{\xi}\circ P_i = P_i\circ A_{\xi},$$
$i = 1,2,3$.  Here $P_i$ is the restriction of $J_i$ to $TM$, where $J_i, i=1,2,3$ denotes a local section of the
quaternionic-K\"ahler structure. It is easy to see such hypersurfaces are curvature-adapted. Then a long calculation shows
that the only such hypersurfaces are horospheres and tubes over totally geodesic $\mathbb{H}H^k, 0\leq k < n$. Just as in
the case of $\mathbb{C}H^n$ one may instead apply the analogous spectral data contained in \cite{berndt}, Theorem 4.18  to
avoid this calculation and shorten their proof.

\section{Generalized pseudo-Einstein hypersurfaces}

We come now to the complete classification of generalized pseudo-Einstein hypersurfaces  of $\mathbb{C}H^n$ (Theorem \ref{cor2}).

\proof From the Gauss Equation one calculates that for a real
hypersurface $M\subset \mathbb{C}H^n$, one has
$$
SX = -(2n+1)X  + 3\langle X, U\rangle U +
\mathfrak{m}A_{\xi}(X) - A_{\xi}^2(X),
$$
for $X\in TM$, where $\mathfrak{m} = tr(A_{\xi})$ denotes the mean curvature of $M$.

As has been outlined in the introduction, the only remaining case is that where $M\subset\mathbb{C}H^n$ satisfies the
assumptions of the theorem with $n\geq 3$. The proof of Proposition 5.2, 5.3 and 5.4 in \cite{cecilryan} goes through
(adjusting for the sign of the curvature tensor in the non-compact case), so $M$ is a Hopf hypersurface with at most three
principal curvatures. By B\"oning's theorem we may assume $M$ is degenerate. Then by the spectral data given in Thereom
(\ref{key}) the principal curvatures are $2$ with eigenspace $U$, $1$ with eigenspace $T_{1}$, and $\lambda$ with eigenspace
$T_{\lambda}$ and moreover $JT_{\lambda}\subset T_{1}$. Set $k= dim(T_{\lambda})$: this must be locally constant. It is
standard theory \cite{nieryan} to show the restriction of $A_{\xi}$ to $\mathfrak{D}$ has eigenvalues given as the root of
the equation
$$
A_{\xi}^2 - \mathfrak{m}A_{\xi} + (2n+1 + \rho)Id = 0.
$$
Hence
$$
\mathfrak{m} = 2+ k(\lambda) + (2n- 2- k) = \lambda + 1.
$$
Since $n\geq 3$, $k> 1$ and so one obtains that $\lambda$ is locally
constant. Hence $M$ has constant principal curvatures, and we are
done. \qed

P. Ryan has a proof of this result \cite{pryan} which is independent of B\"oning's work. This involves a detailed
calculation of the Ricci tensor of $\mathbb{C}H^n$ and an investigation of the various possibilities for the dimensions of
$T_{\lambda}$ and $T_{1}$.

\bigskip

\textbf{Acknowledgements} This work was completed as part of a Ph.D. under the supervision of Professor J\"{u}rgen Berndt at
University College Cork, Ireland. The author would like to thank him for his advice, Pat Ryan
for stimulating discussions, and especially the referee, whose comments and suggestions greatly improved this paper. This work was supported by a postgraduate fellowship awarded by the Irish  Research Council for
Science, Engineering and Technology.

\end{document}